\documentclass[reqno,makeidx,12pt]{amsart}
\usepackage{amssymb,amsfonts,amsmath,graphicx
}

\def\C{{\mathbb{C}}}

\def\E{{\mathbb{E}}}
\def\F{{\mathbb{F}}}

\def\N{{\mathbb{N}}}

\def\R{{\mathbb{R}}}

\def\n{{\mathbb{N}}}

\newcommand{\dst}{\displaystyle}

\def\cB{{\mathcal B}}

\def\cI{{\mathcal I}}

\def\intH{\buildrel o \over H}

\def\ds{\displaystyle}
\def\fig{ \centerline{Fig. \the\count200\global\advance\count200 by 1}}
\def\somtrois#1#2#3{\displaystyle 
\mathop{\bigoplus}\limits_{\mathop{#1}\limits_{{\scriptstyle 
\mathop{#2}\limits_{\scriptstyle #3}}}}}

\newtheorem{thm}{Theorem}[section]
\newtheorem{cor}[thm]{Corollary}
\newtheorem{lem}[thm]{Lemma}

\newtheorem*{thm*}{Theorem}
\newtheorem*{lem*}{Lemma}
\newtheorem*{cor*}{Corollary}

\theoremstyle{definition}
\newtheorem{defn}[thm]{Definition}

\newtheorem{exs}[thm]{Examples}

\theoremstyle{remark}
\newtheorem{rem}[thm]{Remark}

\numberwithin{equation}{section}

\def\croi{\hbox{$\mathop{\mathrel\times\joinrel\mathrel{\vrule height 5pt
 depth 0pt}}$}}

\title[Ergodic theorems for free group actions]{On  ergodic theorems for free group actions
on noncommutative spaces}

\author{Claire Anantharaman-Delaroche}
\address{D\'epartment de Math\'ematiques, Universit\'e d'Orl\'eans,
B. P. 6759, F-45067 Orl\'eans Cedex 2}
\email{claire@labomath.univ-orleans.fr}
\dedicatory{} 

\subjclass{Primary 46L53, 46L55; Secondary 46L50}

\keywords{Noncommutative ergodic theorems, Free group actions}
\date{}

\begin{document}
\maketitle

\begin{abstract}
 We extend in a noncommutative setting the individual ergodic theorem of Nevo and Stein concerning
 measure preserving actions of free groups and averages on spheres $s_{2n}$ of even radius.
 Here we study state preserving actions of free groups on a von Neumann algebra $A$
 and the behaviour of $(s_{2n}(x))$ for $x$ in  noncommutative spaces $L^p(A)$. For the Ces\`aro means
 $\frac{1}{n}\sum_{k=0}^{n-1} s_k$ and $p = +\infty$, this problem was solved by Walker.
 Our approach is based on ideas of Bufetov. We prove a noncommutative version of Rota
 ``Alternierende Verfahren'' theorem. To this end, we introduce specific dilations of the
 powers of some noncommutative Markov operators.
\end{abstract}

\section{Introduction}
Ergodic theorems for measure preserving group actions have a long history,
going back to Birkhoff and von Neumann around 1930. Until some ten years ago they mainly concerned
actions of amenable groups. Initiated by Arnold and Krylov \cite{AK}   and 
carried on by Guivarc'h \cite{Gui}, Grigorchuk \cite{Gri} and Nevo \cite{Nev}
  for free group actions, the ergodic theory of actions of non amenable groups 
has greatly progressed since,
due to the works of Nevo, Nevo-Stein, and Margulis-Nevo-Stein. In particular, these authors
have developed powerful techniques well adapted to actions of semi-simple groups.

It is interesting, and not straightforward, to study the analogues of these 
results in noncommutative (or quantum) probability theory. The setting is a 
($W^\star$-)noncommutative probability space, that is, in this paper,
a pair $(A,\varphi)$ where $A$ is a von Neumann algebra and $\varphi$ a faithful
normal state on $A$. The first result obtained in this framework
is the noncommutative individual ergodic theorem due to Lance \cite{Lan}. 
He proved that, for any automorphism $\sigma$
of $A$ which preserves $\varphi$,  the averages
$$c_n(x) = \frac{1}{n}\sum_{k=0}^{n-1} \sigma^k(x)$$
converge almost uniformly to a $\sigma$-invariant element in $A$. The almost uniform convergence, which will be
defined below (see definition \ref{alsure}),  is the noncommutative analogue of almost everywhere convergence. 
Shortly after, the case of
amenable group actions in the noncommutative framework was studied by Conze and Dang Ngoc \cite{CDN}.
Lance's result was further
extended by K\"{u}mmerer \cite{Kum} (replacing $\sigma$ by any normal positive map
such that $\sigma(1) \leq 1$ and $\varphi\circ \sigma\leq \varphi$) and 
Yeadon \cite{Yea} who studied the convergence of $\big(c_n(x)\big)$ when $x$
belongs to the $L^1$-space $L^1(A,\tau)$, in case $\tau$ is a 
faithful normal semi-finite trace. For more
informations on the state of subject before 1990  we refer to the above
mentioned papers and to the books of Jajte \cite{Jaj}, \cite{Jaj1}.
 The most recent and important developments
in that direction are due to Junge and Xu \cite{JX, JX3}. As a consequence of their maximal ergodic
inequalities in noncommutative $L^p$-spaces,
they generalized in particular Yeadon's theorem to every space $L^p(A,\tau)$, for $p\geq 1$.   

The case of free group actions in quantum probability was considered
by Walker in \cite{Wal}. In order to state his result, as well as the results
of Nevo \cite{Nev} and Nevo-Stein \cite{NS} relative to measure preserving free group actions
we need to introduce some notations. We shall denote by $\F_d$ the free group with $d$ generators $g_1,\dots
,g_d$ and
by $|w|$ the length of $w \in \F_d$, that is the smallest number of generators, together with their inverses,
needed to write the word $w$. Given $d$ automorphisms $\sigma_1, \dots, \sigma_d$ of a von Neumann
algebra $A$, we consider the corresponding homomorphism $w \mapsto \sigma_w$ from $\F_d$ into the group
$\text{Aut}(A)$ of automorphisms of $A$ defined by assigning $\sigma_i$ to $g_i$, $1\leq i\leq d$.
The averaging  operator $s_n$ on the sphere $\mathcal{S}_n$
of radius $n$ is defined as
$$s_n(x) = \frac{1}{\#\mathcal{S}_n}\sum_{w\in \mathcal{S}_n} \sigma_w(x),$$
for $x\in A$,
where $\mathcal{S}_n = \{w\in \F_d : |w| = n\}$.
Let us recall the Nevo-Stein theorem. 

\begin{thm*}[\textbf{Nevo-Stein, \cite{Nev,NS}}]
Let $A = L^\infty(X,\mu)$ be the commutative von Neumann algebra associated with 
a probability space $(X,\mu)$. Let $\sigma_i$, $1\leq i \leq d$, be automorphisms of $A$
induced by measure preser\-ving transformations of $X$.  Let $f\in L^p(X,\mu)$. Then
\begin{itemize}
\item[(i)] If $p>1$, the sequence $\big(s_{2n}(f)\big)$ 
converges almost everywhere and in $L^p(X,\mu)$ to the conditional
expectation $\E(f|\cI_2)$  with respect to the $\sigma$-field $\cI_2$ 
of $\F_{d}^{(2)}$-invariant measurable subsets, where $\F_{d}^{(2)}$
is the subgroup of even length words.
\item[(ii)]  If $p\geq 1$, the sequence $\frac{1}{n} \sum_{k=0}^{n-1} s_k(f)$
converges almost everywhere and in $L^p(X,\mu)$ to the conditional
expectation $\E(f|\cI)$  with respect to the $\sigma$-field $\cI$ 
of $\F_{d}$-invariant measurable subsets.
\end{itemize}
\end{thm*}  

In the noncommutative setting, Walker obtained the following ge\-ne\-raliza\-tion of a part
of the Nevo-Stein theorem.

\begin{thm*}[\textbf{Walker, \cite{Wal}}] Let $\sigma_i$, $1\leq i\leq d$, be automorphisms of a
von Neumann algebra $A$ which leave a faithful normal state $\varphi$ invariant.
For $x\in A$, the sequence  $\frac{1}{n} \sum_{k=0}^{n-1} s_k(x)$ converges almost
uniformly to an element $\hat{x}\in A$.
\end{thm*}

The proof given by Walker uses the result of Nevo and Stein showing that the Ces\`aro means
$\frac{1}{n} \sum_{k=0}^{n-1} s_k$ are dominated by Ces\`aro averages of powers
of the only contraction $s_1$. In the commutative case, Nevo and Stein could conclude
with the help of the classical Hopf-Dunford-Schwartz maximal inequality. 
In the noncommutative case, this inequality has to be replaced by a maximal inequality
due to Goldstein. In both cases, delicate spectral estimations are also needed.

Recently, Bufetov \cite{Buf} has proposed another proof of the Nevo-Stein theorem. In addition to
being very simple, another advantage of his method
is that it allows to extend part (i) of the Nevo-Stein theorem to functions $f$
in the class $L\log L$.

The aim of this paper is to show how the method of Bufetov can be adapted
to quantum probability theory. Combined with recent noncommutative
martingale convergence results due to Junge \cite{Jun} and to Defant and Junge \cite{DJ}, it gives
the following result.

\begin{thm*}[\textbf{Noncommutative Nevo-Stein ergodic theorem}] Let $\sigma_i$,\\ $1\leq i\leq d$, be
automorphisms of a von Neumann algebra $A$ which leave a faithful normal state $\varphi$ invariant.
Let $x\in L^p(A)$ with $p\in ]1,+\infty]$.
\begin{itemize}
\item[(i)] the sequence $s_{2n}(x)$ converges bilaterally almost surely
to the conditional expectation of $x$ with respect to the 
 $\F_{d}^{(2)}$-invariant elements in $L^p(A, \varphi)$.
\item[(ii)] the sequence  $\frac{1}{n} \sum_{k=0}^{n-1} s_k(x)$ converges bilaterally almost
surely to the conditional expectation of $x$ with respect to the 
 $\F_{d}$-invariant elements in $L^p(A, \varphi)$.
\end{itemize}
Moreover, in both cases, the convergence holds almost surely whenever $p\in [2,+\infty[$
and almost uniformly when $p=+\infty$.
\end{thm*}
 
The various notions of convergence appearing in this statement are analogues of almost
everywhere convergence. They are defined in \ref{alsure}.

 This paper is organized as follows.
In the next section, we shall begin by explaining  the notions of noncommutative probability theory
 used to state and prove the above theorem. 
 The reader only interested in a very short proof of the result of Walker may skip
 these technical preliminaries and pass directly to section 3. In theorem \ref{cesaro} we show 
 that the techniques of Bufetov developed in \cite{Buf1} apply immediately to extend the theorem of Walker. The idea
 is to relate the problem to the study of an appropriate Markov operator $P : A^{2d} \to A^{2d}$
 and to use available ergodic theorems for the sequence $\frac{1}{n}
 \sum_{k=0}^{n-1} P^k$ of averaging operators .
 
Part (i) of the noncommutative Nevo-Stein theorem requires a more involved
 analysis. In section 4 we first construct a suitable ``noncommutative Markov chain''
 attached to the Markov operator $P$ (see theorem \ref{basic}) and we prove (corollary \ref{almunif})
 a noncommutative version of Rota ``Alternierende Verfahren'' theorem, that is a convergence
 result for the sequence $P^n \circ (P^\star)^n$ where $P^\star$ is an adjoint of $P$ (as defined in section 2).
 Following Bufetov, we show in  section 5 how this last result implies
 the noncommutative Nevo-Stein theorem.

The above operator $P$ is factorizable in the sense of definition \ref{factor}. Roughly speaking,
 this means that $P = j_{0}^\star\circ j_1$, where $j_0, j_1$ are homomorphims preserving
 given states in a strong sense (see definition \ref{station}). Note that in the commutative case,
any measure preserving Markov operator is factorizable.
In the last section,  we show how a Daniell-Kolmogorov type extension
of   the classical construction
of the Markov chain associated with a transition probability and an initial distribution
can be carried out  for every factorizable Markov operator.
 Such constructions  already appeared in the theory of quantum stochastic processes
 (see \cite{Sau2} and \cite{BP} for instance),
 but here we insist on having, in addition, a good behaviour with respect to $P^\star$
 (see formula (6.3)).
In particular, the noncommutative Rota theorem holds for such an operator $P$.
However, since the construction involves amalgamed free products of
von Neumann algebras, we have chosen to provide in section 4 a much simpler construction
for the concrete Markov operator used in the proof of the noncommutative Nevo-Stein theorem.

I am grateful to Marius Junge and Jean-Luc Sauvageot for useful discussions.
I am specially indebted to Marius Junge for having communicated to me preliminaries
versions of \cite{DJ}, \cite{JX3} and \cite{JX1}. In particular, the proof of lemma \ref{convergence}
in the noncommutative case is due to him.

\section{Preliminaries}
 
 Let us briefly summarize the main concepts and results needed in this paper. 

\subsection{ Tomita-Takesaki modular theory and noncommutative $L^P$-spaces}
We refer to
 \cite{Tak}, \cite{KR}, \cite[Section 2.5]{BR} for general backgrounds on the modular theory 
of von Neumann algebras and to \cite{Haa}, \cite{Kos}, \cite{Ter} 
for details on noncommutative $L^p$-spaces.

When $\tau$ is a normal faithful semi-finite trace on a von Neumann algebra $A$, the
spaces $L^p(A,\tau)$ are easy to introduce and well understood (see \cite{Seg}, and \cite{Nel}
for a short exposition). 
 
 Let us concentrate on the case of a noncommutative probability space $(A,\varphi)$ where $\varphi$
is a normal faithful state. We denote by
$L^2(A,\varphi)$  the completion of
 $A$ with respect to the scalar product $\langle a_1, a_2\rangle = \varphi(a_{2}^* a_1)$
and by $\xi_\varphi$ the unit of $A$ viewed as a vector of $L^2(A,\varphi)$.
We identify $A$ with its $GNS$-representation in $L^2(A,\varphi)$. In particular,
the inclusion $A\subset L^2(A,\varphi)$ is given by the map $a \mapsto  a\xi_\varphi$.
We denote by $S_\varphi$ the closure of the operator $a\xi_\varphi \mapsto a^*\xi_\varphi$.
As usual, its polar decomposition is written $S_\varphi = J_\varphi \Delta_\varphi$.
Recall that $J_\varphi$ is an anti-unitary involution and that $\Delta_{\varphi}^{it} A \Delta_{\varphi}^{-it}
= A$ for all $t\in \R$. The {\it modular automorphism group} associated with $(A,\varphi)$
is the one-parameter automorphism group of $A$ defined by 
$\sigma_{t}^\varphi (a) = \Delta_{\varphi}^{it} a \Delta_{\varphi}^{-it}$ for $a\in A$ and $t\in \R$.

 There are several ways to introduce noncommutative $L^p$-spaces.
We follow the construction of Haagerup \cite{Haa}. For simplicity, let us set $\sigma_t = \sigma_{t}^\varphi$
for $t\in \R$. Given a concrete representation of $A$ on a Hilbert space $H$
(for instance the GNS-representation in $L^2(A,\varphi)$), 
recall that the crossed product
$A\croi_{\sigma} \R$ is the von Neumann algebra of operators acting on $L^2(\R,H)$
 generated by $\pi(a)$, $a\in A$, and $\lambda(s)$, $s\in \R$, where for $\xi\in L^2(\R,H)$ and
 $t\in \R$,
 $$\pi(a)(\xi)(t) = \sigma_{-t}(a)\xi(t) \quad \text{and} \quad \lambda(s)(\xi)(t) = \xi(t-s).$$
 This crossed product does not depend on the choice of $H$. Since $\pi$ is a normal
 faithful representation of $A$ on $L^2(\R,H)$, one identifies
 $A$ with $\pi(A)$. There is a one parameter automorphism group $ t\mapsto \hat{\sigma}_t$
  of $A\croi_{\sigma} \R$, implemented by the unitary representation $t\mapsto W(t)$ of
 $\R$ on $L^2(\R,H)$, where $W(t)(\xi)(s) = e^{-its}\xi(s)$ for $\xi \in L^2(\R,H)$ and
 $s,t \in \R$. One has 
 $$
 \hat{\sigma}_t(\lambda(s)) = W(t) \lambda(s) W(t)^\star = e^{-ist}\lambda(s) \quad\text{for} \quad s,t\in \R,
 $$
 and
 \begin{equation}\label{fix}
 A = \{ x\in A\croi_{\sigma} \R : \hat{\sigma}_t(x) = x, \forall t\in \R\}.
 \end{equation}
 
 The crossed product $A\croi_{\sigma} \R$ has a canonical normal semi-finite faithful trace
 $\tau$ satisfying $\tau\circ \hat{\sigma_t} = e^{-t} \tau$ for all $t\in \R$. We can therefore introduce
 the topological $*$-algebra $\mathcal{M}(A\croi_{\sigma} \R,\tau)$ formed of  the closed densely defined
 operators on $L^2(\R,H)$, affiliated with $A\croi_{\sigma} \R$, that are measurable 
 with respect to $\tau$ \cite[Chapter IX]{Tak}. This algebra is the substitute for the algebra of measurable functions
 in the commutative case. Following the point of view of Haagerup, for $p\in [1,+\infty]$
 we define the noncommutative $L^p$-space $L^p(A,\varphi)$ as
 $$L^p(A,\varphi) = \{x\in \mathcal{M}(A\croi_{\sigma} \R,\tau) : \hat{\sigma}_t(x) = e^{-t/p}x, \forall t\in \R\}.$$
 It is an ordered Banach  space.  We shall not describe its norm $\|.\|_p$ here (see \cite{Haa}).
 Its positive cone $L^p(A,\varphi)_+$ is the intersection of $L^p(A,\varphi)$ with the cone formed
 by the positive mesurable operators. Since $L^p(A,\varphi)$ does not depends on $\varphi$, 
 up to order preserving isometry, we shall write $L^p(A)$ for its  abstract version. These $L^p$-spaces
 behave as their commutative analogue with respect to duality.

Note that we have given above two definitions of $L^2(A,\varphi)$. This is not confusing since,
due to the unicity of the standard form of a von Neumann algebra, the two spaces can be
identified in a natural way.
 
 Observe also  that, as a consequence of (\ref{fix}), we have $L^\infty(A,\varphi) = A$. On the other
 hand, $L^1(A,\varphi)$ is identified we the predual $A_*$ of $A$ in the following way:
 any normal positive linear form $\omega$ on $A$ induces a dual weight $\hat \omega$ on $A\croi_{\sigma} \R$
 whose Radon-Nikodym derivative $h_\omega$ with respect to $\tau$ belongs to $L^1(A,\varphi)$
 and this map $\omega \mapsto h_\omega$ extends to an order preserving isometry 
 between $A_*$ and $L^1(A,\varphi)$.
 
 As concrete subspaces of $\mathcal{M}(A\croi_{\sigma} \R,\tau)$, the 
 $L^p$-spaces satisfy the odd relation
 $$L^p(A,\varphi) \cap L^q(A,\varphi) = \{0\} \quad \text{if} \quad p\not= q.$$
 However, one can define nice positive embeddings as follows. Let us denote by $D = D_\varphi$ ($=h_\varphi$)
 the Radon-Nikodym derivative  of the dual weight $\hat \varphi$ with respect
 to $\tau$.
 
 \begin{lem}[Theorem 1.7, \cite{GL}]\label{embed} Let $p\in [1,+\infty]$ and let $p'$ be such $1/p + 1/p' = 1$.
 \begin{itemize}
 \item[(i)] $\iota_p : x\mapsto D^{\frac{1}{2p}} x D^{\frac{1}{2p}}$ is an embedding from
 $A_+$ into $L^p(A,\varphi)_+$ with dense range;
 \item[(ii)] $\kappa_p : x\mapsto D^{\frac{1}{2p'}} x D^{\frac{1}{2p'}}$ is an embedding from
 $L^p(A,\varphi)_+$ into $L^1(A,\varphi)_+$ with dense range.
 \end{itemize}
 \end{lem}

  \subsection{Markov operators in noncommutative probability theory}
In the commutative case, there are several variant of the notion of Markov operator. We shall adopt a definition
well suited to be extended  to the noncommutative setting.

\begin{defn}[\cite{Kai}]\label{Markov}  A {\it Markov operator} on a probability space $(X,\mu)$ is a positive unital normal
operator $Q$ from $L^\infty(X,\mu)$ into itself. We say that $\mu$ is $Q$-{\it stationary}, or that $Q$
is $\mu$-{\it preserving} if $\ds \int_X Q(f)d\mu = \int_X f d\mu$ for every $f\in L^\infty(X,\mu)$.
\end{defn}

\begin{rem}\label{prolonge1} Usually, Markov operators are defined to be positive contractions from $L^1(X,\mu)$ into itself,
preserving the constant function $1$ (see \cite[page 178]{Neveu}). For $\mu$-preserving operators
the two definitions coincide. Indeed, let $Q$ be a $\mu$-preserving operator in the sense of definition
\ref{Markov}. Then the (predual) operator $Q_* : L^1(X,\mu) \to L^1(X,\mu)$ satisfies
$Q_*(1) = 1$ and therefore preserves the
subspace $L^\infty(X,\mu)$. Hence, the dual of $Q_{*|_{L^\infty(X,\mu)}}$ gives a unique
extension of $Q$ to a positive contraction of $L^1(X,\mu)$.
\end{rem}

A normal unital completely positive map from a noncommutative probability space $(A, \varphi)$ into
 another one $(B,\psi)$ will usually be called a {\it Markov operator}.
As in the commutative case, $Q$ can be extended to  $L^p$-spaces  if $\psi\circ Q = \varphi$.
More generally we have the following result:

  \begin{lem}[\cite{GL},\cite{JX1}] \label{prolonge} Let $Q$ be 
a normal  positive map from $(A,\varphi)$ into  $(B,\psi)$  with $Q(1)\leq 1$ and $\psi\circ Q \leq \varphi$.
The map $Q_{(p)} : D_{\varphi}^{\frac{1}{2p}} A
D_{\varphi}^{\frac{1}{2p}}
\to
  D_{\psi}^{\frac{1}{2p}}BD_{\psi}^{\frac{1}{2p}}$, defined by
  $$Q_{(p)}(D_{\varphi}^{\frac{1}{2p}} a D_{\varphi}^{\frac{1}{2p}}) = D_{\psi}^{\frac{1}{2p}}Q(a)D_{\psi}^{\frac{1}{2p}}$$
  for $a\in A$, extends to a positive contraction $Q_{(p)}$ from $L^p(A,\varphi)$ onto $L^p(B, \psi)$.
  \end{lem}

Hereafter, we shall drop the subscript $p$ and therefore write $Q$ instead of $Q_{(p)}$.
  
  The case of conditional expectations will be especially useful hereafter.
Given $(A,\varphi)$ as above, we shall say that  a von Neumann subalgebra $A_1$ of $A$ is 
$\varphi$-{\it invariant} if it is invariant under the modular
automorphism group $\sigma_{t}^\varphi$ of $\varphi$ ({\it i. e.} 
$\sigma_{t}^\varphi(A_1) \subset A_1$ for every $t\in \R$). 
This condition is equivalent (see \cite{Tak1}) to the existence of a (unique) normal conditional
expectation $\E : A \to A_1$ such that $\varphi \circ \E = \varphi$. 
Moreover, in this situation, we have $\E\circ \sigma_{t}^\varphi = \sigma_{t}^\varphi \circ \E$ for
$t\in \R$. If $\varphi_1$ denotes the restriction of $\varphi$ to $A_1$,
  the space $L^p(A_1, \varphi_1)$ is naturally embedded as a Banach subspace  of
  $L^p(A,\varphi)$ and $\E$ extends to a contractive projection from $L^p(A,\varphi)$
  onto $L^p(A_1, \varphi_1)$, still denoted $\E$ (see lemma \ref{prolonge} or \cite[Section 2]{JX2}).
For $x\in L^p(A,\varphi)$, we shall say that  $\E(x)$ is the {\it conditional expectation
of $x$ with respect to} $A_1$.

In the commutative case, given a $\mu$-preserving Markov operator $Q$ there is a unique
Markov operator on $(X,\mu)$, that we shall denote $Q^\star$, such that $\ds \int_X Q^\star(f)g d\mu
= \int_X f Q(g) d\mu$ for every $f,g \in L^\infty(X,\mu)$ (take $Q^\star =  Q_{*|_{L^\infty(X,\mu)}}$ in 
remark \ref{prolonge1} above). In general the situation is more subtle, and we have to take the
modular automorphism groups of the noncommutative probability spaces into account. Recall
that for any normal unital completely positive map $Q$ from $(A,\varphi)$  to $(B,\psi)$
 such that $\psi\circ Q = \varphi$,
there always exists a unital completely positive map $Q^\star: B \to A$ with 
\begin{equation}\label{adjoint1}
\varphi\big(Q^\star(b)\sigma_{-i/2}^\varphi(a)\big)=
\psi(\sigma_{i/2}^\psi(b)Q(a))
\end{equation}
 for every $\sigma^\varphi$-analytic element $a\in A$
and every $\sigma^\psi$-analytic element $b\in B$ (see \cite{AC},
\cite{CNT}). There exists $Q^\star$ such that 
 \begin{equation}\label{adjoint}
\varphi\big(Q^\star(b)a\big) = \psi\big(bQ(a)\big)
\end{equation}
 for all $a\in A, b\in B$  if and only if $Q$ intertwines the
 modular automorphism groups of $\varphi$ and $\psi$.
We give below a proof of this result, for the reader's convenience.

 \begin{lem}[Proposition 6.1, \cite{AC}]\label{commute} Let $Q$ be a completely positive normal unital map from
$(A,\varphi)$ into $(B,\psi)$. The two following conditions are
equivalent:
   \begin{itemize}
   \item[(i)] there exists a normal unital completely positive map $R : (B,\psi) \to (A,\varphi)$
   such that $\varphi(R(b)a) = \psi(bQ(a))$ for every $a\in A$ and $b\in B$;
   \item[(ii)]  $\psi\circ Q = \varphi$ and
   $\sigma_{t}^\psi \circ Q = Q \circ \sigma_{t}^\varphi$
for every $t\in \R$.
   \end{itemize}
  \end{lem}
  
  \begin{proof}  Note first that when $Q$ a is unital completely positive map with $\psi\circ Q = \varphi$,
   there exists a unique contraction $V : L^2(A,\varphi)\to L^2(B,\psi)$
   such that $V(a\xi_\varphi) = Q(a)\xi_\psi$ for $a\in A$. Moreover,
   we have $V S_\varphi \subset S_\psi  V$.
  
   Assume first the existence of $R$ as in $(i)$. In particular we have $\psi\circ Q = \varphi$
  and $\varphi \circ R = \psi$.
  It is easily checked that the adjoint $V^*$ of $V$ is the operator constructed similarly from
   $R$.  
   We have   $V^* S_\psi \subset S_\varphi  V^*$, so that $VS_\varphi^{*}S_\varphi \subset
  S_\psi^{*}S_\psi V$ and therefore
 $V \Delta_{\varphi}^{it} = \Delta_{\psi}^{it} V$ for every $t\in \R$.
 We can conclude that $\sigma_{t}^\psi \circ Q = Q \circ \sigma_{t}^\varphi$ since
 $$\sigma_{t}^\psi\circ Q(a)\xi_\psi = \Delta_{\psi}^{it} V(a\xi_\varphi)
 =V  \Delta_{\varphi}^{it}(a\xi_\varphi) = Q \circ \sigma_{t}^\varphi(a)\xi_\psi$$
 for $a\in A$ and $t\in \R$.
  
  Let us now prove  that $(ii) \Rightarrow (i)$. Since $\sigma_{t}^\psi \circ Q = Q \circ \sigma_{t}^\varphi$
  we get $V \Delta_{\varphi}^{it} = \Delta_{\psi}^{it} V$ for every $t\in \R$
  and therefore $J_\psi V = V J_\varphi$. Obviously we have $V(A_+ \xi_\varphi) \subset B_+ \xi_\psi$.
  For $b\in B_+$ and $a'  = J_\varphi a J_\varphi \in A_{+}'$, observe that
  $$
  \langle V^* b\xi_\psi, a' \xi_\varphi \rangle = \langle b\xi_\psi, J_\psi Va\xi_\varphi \rangle \geq 0.$$
  It follows that 
  $V^*(B_+ \xi_\psi) \subset \overline{A_+ \xi_\varphi}$ since
$$\overline{A_+ \xi_\varphi}= \{\eta \in L^2(A,\varphi) : \langle \eta,
  a'\xi_\varphi \rangle \geq 0, \forall a'\in A_{+}'\}$$
  (see \cite[Proposition 2.5.27]{BR}). 
    
   Since $\|V^*\| =1$ and $\langle V^* \xi_\psi, \xi_\varphi \rangle = 1$, we see that 
   $V^* \xi_\psi= \xi_\varphi$. It follows from \cite[Lemma 3.2.19]{BR} that
   $V^*(B_+ \xi_\psi) \subset A_+ \xi_\varphi$. 

For $b\in B_+$, let us define
   $R(b)$ as the unique element $a\in A_+$ such that $V^* b\xi_\psi = a \xi_\varphi$.
   Then it is easily checked that $R$ fulfils the conditions of $(i)$.
 \end{proof}
   
    \begin{defn}\label{station}  Let $Q$ be a completely positive normal unital map from
$(A,\varphi)$ into $(B,\psi)$.  If the equivalent conditions of the previous lemma are satisfied, we say
   that the pair $(\psi, \varphi)$ is {\it stationary with respect to $Q$} or  $Q$-{\it stationary}.
   We shall also say that $Q$ is a $(\psi,\varphi)$-{\it preserving Markov operator}.
  The operator $Q^\star$ will be called
the {\it adjoint of $Q$ with respect to $(\psi,\varphi)$}. 

   When $(A,\varphi) = (B,\psi)$, we say that $\varphi$ is {\it stationary with respect to $Q$} or $Q$-{\it stationary},
   or that $Q$ is a {\it $\varphi$-preserving Markov operator}.
\end{defn}

{\sl We insist on the fact that $Q$-stationarity is strictly stronger (in general) than the equality $\psi\circ Q = \varphi$.}

\begin{exs}\label{examples} {\bf (a)} {\it Finite von Neumann algebras.} The modular theory is trivial if $\varphi$ and $\psi$ are 
chosen to be tracial normal faithful states. Therefore, the maps $Q$ that we consider are normal
completely positive unital maps with $\psi\circ Q = \varphi$.

{\bf (b)} {\it Conditional expectations.} As already observed, normal conditional expectations $\E$ with
$\varphi\circ \E = \varphi$ automatically commute with $\sigma^\varphi$. 

{\bf (c)} {\it Homomorphisms.} Let $Q$ be a unital normal $*$-homomorphism from
$A$ into $B$ such that $\psi\circ Q = \varphi$. Then $Q$ satisfies the equivalent conditions
of lemma \ref{commute} if and only if $Q(A)$ is invariant under the modular automorphism group
of $\psi$. Moreover in this case, for $b\in B$ we have
   $Q^\star(b) = Q^{-1}\big(\E(b)\big)$, where $\E$ is the $\psi$-preserving conditional
   expectation from $B$ onto $Q(A)$.
\end{exs}
 
  \subsection{Almost sure convergence in noncommutative probability theory} 
  Finally, let us introduce noncommutative substitutes for almost sure convergence
  (see \cite{Jaj}, \cite{Jaj1} and \cite{DJ}).
  
\begin{defn}\label{alsure} Let $A$ be a von Neumann algebra with a faithful normal state $\varphi$.
  \begin{itemize}
\item[(a)]  We say that a sequence $(x_n)$ of elements of $A$ converges to $0$
{\it almost uniformly} (resp. {\it bilaterally almost uniformly}) if 
for every $\varepsilon >0$ there is a projection $e\in A$
with $\varphi(1-e)\leq \varepsilon$ and $\lim_{n\to \infty}\|x_ne\|_\infty = 0$
  (resp. $\lim_{n\to \infty}\|ex_n e\|_\infty = 0$).
\item[(b)]  Let $p\in [1,+\infty]$. We say that a sequence $(x_n)$ of elements of $L^p(A,\varphi)$ converges to $0$
  {\it almost surely} if for every $\varepsilon >0$ there is a projection $e\in A$ and a family $(a_{n,k})$ in $A$
 such that 
 $$\varphi(1-e)\leq \varepsilon, \, \,x_n = \sum_k (a_{n,k} D^{1/p}),
 \,\,\hbox{and}  \,\lim_{n\to \infty} \Big\|\sum_k ( a_{n,k}e)\Big\|_\infty = 0,$$
 where the two series converge in $L^p(A,\varphi)$ and $A$ respectively.
  \item[(c)]   Let $p\in [1,+\infty]$. We say that a sequence $(x_n)$ of elements of $L^p(A,\varphi)$ converges to $0$
  {\it bilaterally almost surely} if for every $\varepsilon >0$ there is a projection $e\in A$ and a family $(a_{n,k})$ in $A$
 such that 
 $$\varphi(1-e)\leq \varepsilon, \,\, x_n = \sum_k (D^{1/2p}a_{n,k} D^{1/2p}),
 \,\,\hbox{and} \, \lim_{n\to \infty} \Big\|\sum_k (e a_{n,k}e)\Big\|_\infty = 0,$$
 where the two series converge in $L^p(A,\varphi)$ and $A$ respectively.
\end{itemize}
\end{defn}

When $A = L^\infty(X,m)$, with $(X,m)$ a probability measure,  all these notions
of convergence coincide, via Egorov's theorem, with the almost everywhere pointwise convergence.
When $A$ is a semi-finite von Neumann algebra equipped with a normal faithful semi-finite trace
 $\tau$, the notions of almost uniform convergence and bilateral almost uniform convergence
 can still be introduced as in definition \ref{alsure} (a) for every sequence in $L^p(A,\tau)$ and all $p\in [1,+\infty]$.
 
\section{Ergodic theorems for some Ces\`aro means}

Let us begin by introducing some notations.
Let $\big[p(ij)\big]$ be a  stochastic matrix, whose rows and columns
are indexed by a finite set $I$ and let $\big(p(i)\big)_{i\in I}$ be a stationary distribution, that is a
probability measure on $I$ such that   $p(j) = \sum_{i\in I} p(i) p(ij)$ for 
$j\in I$. We assume that $p(i) > 0$ for all $i\in I$.

We shall usually view elements $w\in I^{n}$ as words of length $n$ in the alphabet $I$
and write $w = w_0w_1\cdots w_{n-1}$ where $w_k$ is the $(k+1)-th$ component of $w$.
For $w\in I^{n}$,  $n\geq 2$, we set 
$$p_{n-1}(w) = p(w_0)p(w_0w_1) \cdots p(w_{n-2}w_{n-1})$$
 and 
$$I(n) = \{w\in I^{n} : p_{n-1}(w) \not= 0\}.$$
We set $p_0(i) = p(i)$ and $I(1) = I$.

Let $L$ be a linear space and consider, for $i\in I$, linear operators $P_i : L \to L$.
For $n\geq 1$ and $w\in I(n)$, we set $P_w = P_{w_0}\circ\cdots\circ P_{w_{n-1}}$. We define
the operators $s_n$ and $c_n$ by the formulas
$$\aligned
s_{n} & = \sum_{w\in I(n)} p_{n-1}(w) P_w \quad\hbox{if} \quad n\geq 1, \quad s_0 = \hbox{Id}_A,\\
c_n & = \frac{1}{n}\sum_{k=0}^{n-1} s_k  \quad\hbox{if} \quad n\geq 1.
\endaligned$$

For us $L$ will be a noncommutative $L^p$-space and we shall study the convergence
of the sequences $\big(s_n(x)\big)$ and $\big(c_n(x)\big)$, where $x\in L$. 

Following the idea of Bufetov \cite{Buf1}, we assign to our data the operator $P : L^I \to L^I$
such that 
$$P(x)_i = \sum_{j\in I} p(ij) P_i(x_j), \quad i\in I,$$
for all $x = (x_i)_{i\in I}$.

Given $x\in L$, we denote by $\tilde{x}$ the element of $L^I$ with $\tilde{x}_i = x$ for all $i$.
The main observation, easily proved by induction, is that, for $x\in L$ and $n\geq 1$,
\begin{equation}\label{bufetov}
P^n(\tilde{x})_i = \frac{1}{p_i}\sum_{{w\in I(n)}\atop{w_0 = i}} p_{n-1}(w) P_w(x),
\end{equation}
so that
\begin{equation}\label{bufetov1}
s_n(x) = \sum_{i\in I} p(i) P^n(\tilde{x})_i
\end{equation}
and
\begin{equation}\label{bufetov2}
c_n(x) = \sum_{i\in I} p(i) \Big(\frac{1}{n}\sum_{k= 0}^{n-1}P^k(\tilde{x})_i\Big).
\end{equation}

In this section,
we shall concentrate on $\big(c_n(x)\big)$ which is much easier to deal with.

\begin{thm}\label{cesaro} Let $[p(ij)]$ and $\big(p(i)\big)_{i\in I}$ be as above and $(A, \varphi)$ be
a noncommutative probability
space. Let $P_i : A \to A$ be normal positive contractions such that $\varphi \circ P_i \leq \varphi$ ($i\in I$). For 
every $x\in A$, the sequence $\ds c_n(x)  = \frac{1}{n}\sum_{k=0}^{n-1}\Big(\sum_{w\in I(n)} p_{n-1}(w) P_w(x)\Big)$
converges almost uniformly to an element $\hat{x} \in A$.
\end{thm}

\begin{proof} We use the notations introduced above where $L$ is here taken to be $A$.
The space $A^I$ is a von Neumann algebra, that we denote by $B$,  and we equip $B$
with the normal faithful state $\phi$ defined by
$$\phi(b) = \sum_{i\in I} p(i)\varphi(b_i)$$ 
where $b = (b_i)_{i\in I}$.
 
Then $P : B \to B$ defined as above is a normal positive contraction such that $\phi\circ P \leq \phi$.
The ergodic theorem of K\"{u}mmerer \cite{Kum} implies that for every $b\in B$ the sequence
$\ds \Big(\frac{1}{n}\sum_{k=0}^{n-1} P^k(b)\Big)$ converges almost uniformly in $B$
to a $P$-invariant element. It follows from (\ref{bufetov2}) that for every $x\in A$ that sequence $c_n(x)$
converges almost uniformly in $A$.
\end{proof}

\begin{rem}\label{walker1} The situation considered by Walker is the following particular case:
\begin{itemize}
\item[-] $I$ is a set of $2d$ elements, that we write 
 $$I = \{-d, -d+1,\dots,-1,1,\dots, d-1, d\};$$
 \item[-]  for $1\leq i\leq d$, $P_i$ is an automorphism $\sigma_i$ of $A$
 such that $\varphi\circ \sigma_i = \varphi$ and  we set $\sigma_{-i} = \sigma_{i}^{-1}$;
 \item[-] $\big[p(ij)\big]$ is the stochastic matrix  defined by
$$p(ij) =
\begin{cases}
& 0\quad \text{if} \quad i=-j,\\
& \frac{1}{2d-1}\quad \text{otherwise}.
\end{cases}
$$
 \item[-] $\big(p(i)\big)_{i\in I}$ is the uniform probability measure on $I$, {\it i.e.} $p(i) = 1/2d$ for all $i\in I$.
\end{itemize}
  Obviously, in this case, $I(n)$ is the sphere $\mathcal{S}_n$ formed by the words
  of length $n$.
 \end{rem}
 
 The identification of the limit in theorem \ref{cesaro} is a more difficult problem. However, in the situation
 of the previous remark, it is easy to show that the limit $\hat{x}$ is the conditional expectation of $x$
 on the subalgebra of fixed points under the action of $\F_d$ (see the proof of the Nevo-Stein theorem
 in section 4).
 
 Similarly, using the noncommutative individual Dunford-Schwartz ergodic theorem obtained recently by Junge and Xu \cite[Corollaire 5]{JX}, \cite[Corollary 6.4]{JX3}, we get:
 
 \begin{thm} Let $A$ be a semifinite von Neumann algebra equipped with a normal
semi-finite faithful trace $\tau$. Let $P_i : A \to A$, $i\in I$, be normal positive contractions 
such that $\tau\circ P_i \leq \tau$.
 We are given $[p(ij)]$ and $\big(p(i)\big)_{i\in I}$ as in theorem \ref{cesaro}. Then, for every $x\in L^p(A)$ and 
$p\in [1,\infty]$,
 the sequence $\big(c_n(x)\big)$ converges bilaterally almost uniformly to an element of $L^p(A)$.
If $p\geq 2$, the convergence also holds almost uniformly.
 It also holds in $L^p$-norm for $p\in ]1,+\infty[$. 
 \end{thm}
 
 The study of the sequence $(s_n)$ is much more delicate and requires additional assumptions.
 We shall only consider the situation described in remark \ref{walker1}. We need a generalization
 of the Rota ``Alternierende Verfahren'' theorem that we shall examine in the next section.
 
\section{Noncommutative Rota theorem in a simple case}
Let  $(A,\varphi)$ be a noncommutative probability space.
We are given normal  unital endomorphisms $\sigma_i$ of $A$ where  $i$
belongs  to a finite set $I$ of indices. We assume
 that  $\varphi\circ\sigma_i = \varphi$ and  $\sigma_i\circ\sigma_{t}^\varphi =
\sigma_{t}^\varphi\circ \sigma_i$ for $i\in I$ and $t\in \R$. 
We are also given, as before, a  stochastic matrix $\big[p(ij)\big]$
 and a stationary distribution $\big(p(i)\big)_{i\in I}$ with $p(i) > 0$ for all $i\in I$
 
 We set $B= A^I$ and  $\phi$ will be the state defined by
$$\phi(b) = \sum_{i\in I} p(i)\varphi(b_i).$$ 
We still introduce the Markov operator $P$ from $B$ into $B$ 
defined by
\begin{equation} \label{P} 
 P(b)_i = \sum_{j\in I} p(ij)\sigma_i(b_j),
 \end{equation}
for all $b\in B$ and $i\in I$. Each $\sigma_i$ has an adjoint $\sigma_{i}^\star$
with respect to $\varphi$, and therefore $P$ is a $\phi$-preserving Markov operator. We have,
for $b\in B$ and $i\in I$,
$$P^\star(b)_i = \frac{1}{p(i)}\sum_{j\in I} p(j)p(ji) \sigma_{j}^{\star}(b_j).$$

We need the following dilation result.

\begin{thm}[{\bf Noncommutative Kolmogorov construction}] \label{basic}
There exist 
\begin{itemize}
\item [--] a von Neumann algebra $\mathcal{B}$,
 \item [--] a normal faithful state $\Phi$ on $\mathcal{B}$,
  \item [--] a normal unital  endomorphism $\beta : \mathcal{B} \to \mathcal{B}$ with $\Phi\circ \beta = \Phi$
and $\sigma_{t}^\Phi \circ \beta
= \beta \circ \sigma_{t}^\Phi$ for $t\in \R$,
 \item [--] a normal
unital homomorphism $J_0 : B \to \mathcal{B}$ with $\Phi\circ J_0 = \phi$ and $\sigma_{t}^\Phi \circ J_0
= J_0 \circ \sigma_{t}^\phi$ for $t\in \R$,
\end{itemize}
such that, if we set
  $J_n = \beta^n\circ J_0$ for $n\geq 0$, and if $\mathcal{B}_{n]}$ and $\mathcal{B}_{[n}$
denotes the von Neumann subalgebras of $\mathcal{B}$ generated by $\bigcup_{k\leq n}J_k(B)$
 and $\bigcup_{k\geq n}J_k(B)$ respectively, then
 \begin{itemize}
 \item[(i)] the algebras $\mathcal{B}_{n]}$ and $\mathcal{B}_{[n}$ are $\Phi$-invariant;
\item[(ii)] if $\mathbb{E}_{n]}$ and $\mathbb{E}_{[n}$ are the corresponding $\Phi$-preserving
conditional expectations, for $n\in \mathbb{N}$ and $q\geq n$, we have
 \begin{align}\label{basicno}
& \mathbb{E}_{n]}\circ  J_q = J_n \circ P^{q-n},\\
& \mathbb{E}_{[n}\circ J_0 = J_n \circ (P^{\star})^n.
\end{align}
\end{itemize}
 \end{thm}
 
 Let us show first how the following corollary can be deduced from this theorem.
 
\begin{cor}[{\bf Noncommutative Rota theorem}]\label{almunif} For $p\in ]1,+\infty]$ and $x\in L^p(B)$, the sequence
$\big (P^n\circ (P^\star)^n(x)\big )$ converges bilaterally almost surely. For $p\in [2,+\infty]$
 the convergence also holds almost surely.
\end{cor}
 
 \noindent{\it Proof.} We have  $\E_{0]}\circ \E_{[n} \circ J_0(x) =
J_0 \circ P^n \circ (P^{\star})^{n}(x)$ for $x\in B$ and, using lemma \ref{prolonge},
this formula extends to $L^p(B)$, $p > 1$. Therefore the result is a consequence of the lemma below.

\begin{lem}\label{convergence} Let $(\mathcal{B}, \Phi)$ be a noncommutative probability space and let
$(B_n)$ be a decreasing sequence of $\Phi$-invariant von Neumann
subalgebras. We set $B_\infty = \cap B_n$ and we denote by $\E_n$ the
$\Phi$-preserving conditional expectation from $\mathcal B$ onto $B_n$, $0\leq n \leq \infty$. Let
$Q : \mathcal{B} \to \mathcal{B}$ be a normal completely positive contraction. Then
\begin{itemize}
\item[(a)] for $x\in \mathcal B$ (i.e. case $p= +\infty$), the sequence $\big(Q\circ \E_n(x)\big)$ converges to $Q\circ \E_\infty(x)$ almost
uniformly;
\item[(b)] for $p\in [2,+\infty]$ and $x \in L^p(\mathcal{B},\Phi)$, 
the sequence $\big(Q\circ \E_n(x)\big)$ converges to $Q\circ \E_\infty(x)$ almost surely;
\item[(c)] for $p\in ]1,+\infty]$ and $x \in L^p(\mathcal{B},\Phi)$, 
the sequence $\big(Q\circ \E_n(x)\big)$ converges to $Q\circ \E_\infty(x)$ bilaterally almost surely.
\end{itemize}
\end{lem}

\begin{proof} Replacing $x$ by $x- \E_\infty(x)$, we may assume that $\E_\infty(x) = 0$. 

The convergence of $\big(\E_n(x)\big)$ in the appropriate sense is a result due to
Dang-Ngoc \cite[Theorem 4]{DN} for $p = +\infty$ and to Defant and Junge \cite{DJ} in the two other cases.
The main problem, already appearing in the commutative case, is to show that
$\big(Q\circ\E_n(x)\big)$ still converges in the same sense. The case $p=+\infty$ is immediately
solved because $Q$ is continuous on the norm bounded subsets of $\cB$ when $\cB$ is endowed with the
 topology of almost uniform convergence (see for instance \cite[Proposition 1]{Sau1}).

Let us assume that $p\in ]1,+\infty[$ and let us  first recall the proof in the commutative case.
Thanks to Chebichev inequality, to show that  $\big(Q\circ\E_{n}(x)\big)$ goes to $0$ almost
surely, it is enough to show that 
\begin{equation}\label{cheb}
\lim_m \Big\|\sup_{n\geq m}\big |Q\circ \E_n(x)\big | \Big\|_p = 0.
\end{equation}
But, since $Q$ is a positive contraction, we have
\begin{equation}\label{inegal}
 \Big\|\sup_{n\geq m}\big |Q\circ \E_n(x)\big | \Big\|_p \leq \Big\|Q\big(\sup_{n\geq m} \big |\E_n(x) \big |\big)\Big\|_p \leq 
\Big\|\sup_{n\geq m} \big |\E_n(x) \big |\Big\|_p.
\end{equation}

Therefore, it suffices to prove that $\lim_m\Big\|\sup_{n\geq m} \big |\E_n(x) \big |\Big\|_p = 0$.
The classical Doob maximal inequality gives, for all $m\in \N$,
\begin{equation}\label{doob}
\Big \| \sup_n \big |\E_n\big(\E_{m}(x)\big)\big |Ê\Big\|_p \leq
\frac{p}{p-1}\| \E_{m}(x) \|_p.
\end{equation}
Now, observe that for $n\geq m$, we have
$\E_n \circ \E_{m} = \E_n$, and therefore
$$\Big\| \sup_{n\geq m} \big | \E_n(x)\big |\Big\|_p = 
\Big \| \sup_{n\geq m} \big |\E_n\big(\E_{m}(x)\big)\big |Ê\Big\|_p  \leq \frac{p}{p-1}\| \E_{m}(x) \|_p.$$
Since $\big(\E_m(x)\big)$ goes to $0$ in $L^p$-norm, this concludes the proof in the
commutative case.

When $\cB$ is not commutative, the proof follows the same pattern. However, it is not
a straightforward generalization and several major difficulties arise. 
Given a sequence $(x_n)$ in $L^p(\cB, \Phi)$ a first  problem
is to give a meaning to $\sup_n|x_n|$. To this purpose, Junge has introduced in \cite{Jun} the space 
$L^p(\cB, \ell^\infty)$, derived from Pisier's theory of vector valued noncommutative $L^p$-spaces \cite{Pis}.
It is defined as the space of sequences $(x_n)$ in $L^p(\cB,\Phi)$ such that there exist $a,b \in L^{2p}(\cB,\Phi)$
and $(y_n) \in \ell^\infty(\cB)$ with $x_n = ay_n b$ for all $n$. Equipped with the norm
$$\big\| (x_n) \big \|_{L^p(\cB, \ell^\infty)} = \inf \{\|a\|_{2p} \sup_n \|y_n\|_\infty \| b\|_{2p}\},$$
where the infimum runs over all the possible decompositions, $L^p(\cB, \ell^\infty)$
is a Banach space. As explained for instance in \cite{JX3},
$\big\| (x_n) \big \|_{L^p(\cB, \ell^\infty)}$ can be viewed as a noncommutative analogue of 
$\big\| \sup_n |x_n|\big \|_p$. The Doob-Junge inequality for noncommutative martingales (see \cite{Jun})
is expressed as the existence
of a constant $c_p$ such that for all $x\in L^p(\cB,\Phi)$,
$$\Big\| \big(\E_n(x)\big) \Big \|_{L^p(\cB, \ell^\infty)} \leq c_p \|x\|_p.$$

We also need the two following facts, proved by Defant and Junge in \cite{DJ} :

-- a sequence $(x_n)$ in $L^p(\cB, \Phi)$ goes to zero bilaterally almost uniformly if and only if 
$\lim_m \big\| (x_n)_{n\geq m}\big\|_{L^p(\cB, \ell^\infty)} = 0$;

-- given $Q$ as in the statement of the lemma, we have, for $m\in \N$, 
$$\Big\| (Q(x_n))_{n\geq m}\Big\|_{L^p(\cB, \ell^\infty)} \leq \big\| (x_n)_{n\geq m}\big\|_{L^p(\cB, \ell^\infty)}.$$

Having at hand the noncommutative analogues of (\ref{inegal}) and (\ref{doob}), the proof of statement (c)
proceeds exactly as in the commutative case. Statement (b) can be obtained similarly. The space
$L^p(\cB, \ell^\infty)$ has to be replaced by the right hand sided space $L^p(\cB, \ell_{c}^\infty)$,
defined as the space of sequences $(x_n)$  in $L^p(\cB,\Phi)$ such that there exist $b \in L^{p}(\cB,\Phi)$
and $(y_n) \in \ell^\infty(\cB)$ with $x_n = y_n b$ for all $n$. This space is endowed  with the norm
$$\big\| (x_n) \big \|_{L^p(\cB, \ell_{c}^\infty)} = \inf \{ \sup_n \|y_n\|_\infty \| b\|_{p}\},$$
where the infimum runs over all the possible decompositions. Now the Doob-Junge one-sided
maximal inequality reads as
$$\Big\| \big(\E_n(x)\big) \Big \|_{L^p(\cB, \ell_{c}^\infty)} \leq \sqrt{c_{\frac{p}{2}}} \|x\|_p.$$
  \end{proof}

\begin{proof}[Proof of theorem \ref{basic}] We keep some notations introduced in the previous section.
As already said, we shall usually view elements $w\in I^{n}$ as words of length $n$ in the alphabet $I$
and write $w = w_0w_1\cdots w_{n-1}$ where $w_k$ is the $(k+1)-th$ component of $w$.
For  $n\geq 1$, we set 
$I(n) = \{w\in I^{n} : p_{n-1}(w) \not= 0\}$, with $p_n$ defined in section 3.

We denote by $B_n$ the product von Neumann algebra formed by the maps $b : w\mapsto b_w$
from $I(n+1)$ into $A$.  This von Neumann algebra will be
 equipped with the normal faithful state $\phi_n$ defined by
$$\phi_n(b) = \sum_{w\in I(n+1)} p_n(w) \varphi(b_w).$$
Its modular automorphism group $t\mapsto\sigma_{t}^{\phi_n}$ satisfies
$$\sigma_{t}^{\phi_n}(b)_w = \sigma_{t}^{\varphi}(b_w)$$
for $b\in B_n$ and $w\in I(n+1)$.

Note that $B_0 = B$ and $\phi_0 = \phi$.
We shall need the two following unital injective homomorphisms $\alpha_n$ and $\beta_n$
from $B_{n-1}$ into $B_n$, defined, for $n\geq 1$ and $b\in B_{n-1}$ by
\begin{align}
\alpha_n(b)_w &= b_{w_0\cdots w_{n-1}}\\
\beta_n(b)_w &= \sigma_{w_0}(b_{w_1\cdots w_{n}}).
\end{align} 
  They satisfy the following relations~:
$$\aligned
\alpha_{n+1}\circ \beta_n &= \beta_{n+1}\circ \alpha_n\\
\phi_n\circ \alpha_n & = \phi_{n-1} = \phi_n\circ \beta_n.
\endaligned$$
Obviously, $\alpha_n$ and $\beta_n$ intertwine the modular automorphism groups:
 \begin{equation}\label{commutation}
  \alpha_n\circ \sigma_{t}^{\phi_{n-1}} =
  \sigma_{t}^{\phi_n}\circ \alpha_n\quad  \text{and} \quad \beta_n\circ \sigma_{t}^{\phi_{n-1}} =
  \sigma_{t}^{\phi_n}\circ \beta_n
  \end{equation}
   for all $t\in \R$.

  It follows that the algebras $\alpha_n(B_{n-1})$ and $\beta_n(B_{n-1})$ are $\phi_n$-invariant.

We denote by $(\widetilde{\cB}, \widetilde{\Phi})$ the inductive limit of $(B_n, \phi_n)$,
with respect to the morphisms $\alpha_n : (B_{n-1}, \phi_{n-1}) \to (B_n, \phi_n)$,
and by $\lambda_n$ the canonical injective homomorphism from $B_n$ into $\tilde{\cB}$. 
 By definition, we have
 \begin{equation}\label{alpha}
\lambda_{n+1} \circ \alpha_{n+1} = \lambda_n,
 \end{equation}
 and, using the first equality in (\ref{commutation}), we see that $\lambda_n$ is $(\widetilde{\Phi}, \phi_n)$-preserving.
 
Let  $\beta$ the injective endomorphism of $\widetilde{\cB}$ characterized  by
\begin{equation}\label{beta}
\lambda_{n+1} \circ \beta_{n+1} = \beta\circ \lambda_n
\end{equation}
 for every $n\geq 0$.  Note  that $\beta$  is $\widetilde{\Phi}$-preserving. 
 
Next, we set $J_0 = \lambda_0$ and $J_n = \beta^n \circ J_0$ for $n\geq 0$.
Then $\cB$ is defined to be the von Neumann subalgebra of $\widetilde{\cB}$ generated by $\bigcup_{n\geq 0}
J_n(B)$ and we equip this algebra with the restriction $\Phi$ of the state $\widetilde{\Phi}$.  Finally
we denote by $\cB_{n]}$ and $\cB_{[n}$ the von Neumann subalgebras generated by
$\bigcup_{k\leq n} J_k(B)$ and $\bigcup_{k\geq n} J_k(B)$ respectively.
   Obviously, we have $\beta(\cB) \subset \cB$, and we shall still denote by
   $\beta$ the restriction endomorphism $\beta_{|\cB} : \cB \to \cB$. Since
  $\sigma_{t}^{\Phi} \circ J_n = J_n \circ \sigma_{t}^{\phi}$
   for all $n\in \N$, we see that the algebras $J_n(B)$ are $\Phi$-invariant,
   as well as the algebras $\mathcal{B}_{n]}$ and $\mathcal{B}_{[n}$.

Computations hereafter are straightforward. They use the following observation:
  $J_{q}(b)$, when viewed, for $r\geq q$, as an element of $B_{r]} \subset \lambda_r(B_r)$
  (and therefore, after identification, as an element of $B_r$) is written
  $$J_q(b)= \alpha_r\circ\dots\circ\alpha_{q+1}\circ\beta_q\circ\dots\circ\beta_1(b),$$
  hence
  \begin{equation}\label{Jq}
  J_q(b)_{w_0 \cdots w_r} = \sigma_{w_0 \cdots w_{q-1}}(b_{w_q}).
  \end{equation}
  
Let us first prove that $\mathbb{E}_{n]}\circ  J_q = J_n \circ P^{q-n}$ for $q\geq n$. We have to show that
for $c\in B_{n]}$ and $b \in B$, then
$$\Phi\big(cJ_q(b)\big) = \Phi\big(cJ_n\circ P^{q-n}(b)\big).$$
 We may take $c$ of the form
 $$c = J_{k_1}(c^1)J_{k_2}(c^2) \cdots J_{k_\ell}(c^\ell)$$
 where $c^i \in B$ and $0\leq k_i\leq n$ for $1\leq i \leq \ell$. We work in $B_q$. Thanks to
 formula (\ref{Jq}) we get
  $$\aligned
  \Phi\big(c J_q(b)\big) &= \phi_q\big(cJ_q(b)\big)\\
  &=  \sum_{w\in I(q+1)}p_q(w)\varphi\big(C(w)
 \sigma_{w_0 \cdots w_{q-1}}(b_{w_{q}})\big),
  \endaligned
  $$
where we have set $C(w) = \sigma_{w_0\cdots w_{k_1-1}}(c^{1}_{w_{k_1}})
 \cdots \sigma_{w_0\cdots w_{k_\ell-1}}(c^{\ell}_{w_{k_\ell}})$.
  
  On the other hand, we have

\begin{equation}\label{form1} 
  \phi_q\big(cJ_n \circ P^{q-n}(b)\big)
 = \sum_{w\in I(n+1)}p_n(w)\varphi\Big(C(w)
 \sigma_{w_0 \cdots w_{n-1}}\big(P^{q-n}(b)_{w_{n}}\big)\Big),
\end{equation}
  and by iteration of formula (\ref{P}),
 \begin{equation}\label{form2} 
 P^{q-n}(b)_{w_{n}}= \sum_{v\in I(q-n+1) \atop{v_0 =w_n}} p(v_0v_1)\cdots p(v_{q-n-1}v_{q-n})\sigma_{v_0\cdots
v_{q-n-1}}(b_{v_{q-n}}).
 \end{equation}
 
 Writing the above $v$ as $v = w_n \cdots w_q$ and replacing $P^{q-n}(b)_{w_{n}}$
 by its expression (\ref{form2}) in (\ref{form1}), we immediately get that
 $$\phi_q\big(cJ_q(b)\big) = \phi_q\big(cJ_n \circ P^{q-n}(b)\big).$$
    
The proof of the relation $\mathbb{E}_{[n}\circ J_0 = J_n \circ (P^{\star})^n$ is similar.
We have to check that 
 $$\Phi\big(c J_0(b)\big) = \Phi\big(cJ_n \circ (P^\star)^n(b)\big)$$
 for every $c\in \cB_{[n}$ and $b\in B$. It is enough to take  
    $\displaystyle c\in \Big( \bigcup_{n\leq k \leq r} J_k(B)\Big)''$ with $r\geq n$. We may choose
    $c$ of the form 
    $$c = J_{k_1}(c^1) \cdots J_{k_\ell}(c^\ell)$$
 where $c^i\in B$ and $n\leq k_i\leq r$ for $1\leq i \leq \ell$.  Again, we get
    \begin{equation}\label{debut}
\Phi\big(c J_0(b)\big) =\sum_{w\in I(r+1)}p_r(w)\varphi\big(C(w)b_{w_0}\big)
\end{equation}
 and 
  $$ 
  \Phi\big(cJ_n \circ (P^\star)^n(b)\big)
 = \sum_{w\in I(r+1)}p_r(w)\varphi\Big(C(w)
 \sigma_{w_0 \cdots w_{n-1}}\big((P^\star)^n(b)_{w_{n}}\big)\Big),$$
where  $C(w) = \sigma_{w_0\cdots w_{k_1-1}}(c^{1}_{w_{k_1}})
 \cdots \sigma_{w_0\cdots w_{k_\ell-1}}(c^{\ell}_{w_{k_\ell}})$.

Using the invariance of $\varphi$ by the endomorphisms $\sigma_i$, we obtain
 $$ 
  \Phi\big(cJ_n \circ (P^\star)^n(b)\big)
 = \sum_{w \in I(r+1)}p_r(w)\varphi\big(C'(w_{n}\cdots w_r)(P^\star)^n(b)_{w_{n}}\big)$$
 with $C'(w_{n}\cdots w_r) = \sigma_{w_{n}\cdots w_{k_1-1}}(c^{1}_{w_{k_1}})
 \cdots \sigma_{w_{n}\cdots w_{k_\ell-1}}(c^{\ell}_{w_{k_\ell}})$.
 If we begin by summing on the $n$ first letters of $w$, we get that 
 $ \Phi\big(cJ_n \circ (P^\star)^n(b)\big)$ is equal to
\begin{equation}\label{suite}
  \sum_{w_{n}\cdots w_r \in I(r-n+1)} p_{r-n}(w_{n}\cdots w_r)\varphi\big(C'(w_{n}\cdots
w_r)(P^\star)^n(b)_{w_{n}}\big).
 \end{equation}
 
  Furthermore we have
\begin{equation}\label{puissance}
 (P^\star)^n(b)_{w_n} = \frac{1}{p(w_n)}\sum_{v\in I(n+1) \atop{v_n = w_n}}p_n(v) 
\sigma^{\star}_{v_{n-1}}\circ\cdots\circ \sigma_{v_0}^{\star}(b_{v_0}).
\end{equation}
In the expression (\ref{suite}) of  $\Phi\big(cJ_n \circ (P^\star)^n(b)\big)$ we replace $(P^\star)^n(b)_{w_{n}}$ by the
its expression (\ref{puissance}). It follows from the definition of $\sigma_{i}^\star$
that $ \Phi\big(cJ_n \circ (P^\star)^n(b)\big)$
  is equal to
 $$ \sum_{w\in I(r+1)} p_r(w)\varphi\big(\sigma_{w_0\cdots w_{k_1 -1}}(c^{1}_{w_{k_1}})
 \cdots \sigma_{w_0\cdots w_{k_\ell -1}}(c^{\ell}_{w_{k_\ell}})b_{w_0}\big).
  $$
  Comparing with the expression found in (\ref{debut}), this concludes the proof of theorem \ref{basic}
\end{proof}

\begin{rem} Condition (\ref{basicno}) is sometimes called the {\it Markov property} in pro\-bability theory.
In the classical commutative case, we have, furthermore, the following Markov (or covariance)
property (see \cite[Chapter 1]{Rev}) :
\begin{equation}\label{propMarkov}
\forall n,q \in \N, \quad \beta^n\circ \E_{q]} = \E_{n+q]}\circ \beta^n.
\end{equation}

Here also, in the setting of theorem \ref{basic} it can be proved by an easy explicit computation of  $\E_{n]}$
that the Markov property (\ref{propMarkov}) holds.
\end{rem}

\section{Proof of the noncommutative Nevo-Stein theorem}

\begin{proof} We shall apply the previous section in the particular situation
 described in remark \ref{walker1}. Recall that  $$I = \{-d, -d+1,\dots,-1,1,\dots, d-1, d\},$$
 that $\big[p(ij)\big]$
 is the matrix with $p(ij) = 0$ if $i=-j$ and $p(ij) = \ds \frac{1}{2d-1}$ otherwise, and that
  $\big(p(i)\big)_{i\in I}$ is the uniform probability measure on $I$, {\it i.e.} $p(i) = 1/2d$ for all $i\in I$. 
  Furthermore, $\sigma_i$,
 $1\leq i\leq d$, is a $\varphi$-preserving automorphism of $A$ and $\sigma_{-i} = \sigma_{i}^{-1}$.

We  follow the method used by Bufetov in the commutative case \cite{Buf}.
We introduce $B = A^{2d} = \{(b_i)_{i\in I} : b_i \in A\}$, the state $\phi$ on $B$ and the completely
positive map $P$ as in the previous section, that is
  $$P(b)_i = \frac{1}{2d-1} \sum_{ j\in I\atop{j\not= -i}} \sigma_i(b_j).$$

For $x\in A$ we denote by $\tilde x$ the element of $B$ such that $\tilde{x}_i = x$ for $i\in I$.
 Formula (\ref{bufetov}) becomes
   $$P^n(\tilde{x})_i = \frac{1}{(2d-1)^{n-1}} \sum_{ w\in \mathcal{S}_n\atop{w_0= i}} \sigma_w(x)$$
  for $n\geq 1$, so that  
\begin{equation}\label{basrel}
s_n(x) = \frac{1}{\#\mathcal{S}_n} \sum_{ w\in \mathcal{S}_n} \sigma_w(x) =\frac{1}{2d} \sum_{i\in I} P^n(\tilde{x})_i.
\end{equation}

Let us recall that, for $b\in B$, we have
$$P^\star(b)_i = \frac{1}{2d-1} \sum_{ j\in I\atop{j\not= -i}}  \sigma_{-j}(b_j)
  =\frac{1}{2d-1} \sum_{ j\in I\atop{j\not= i}}  \sigma_{j}(b_{-j})$$
since $\sigma_{j}^\star = \sigma_{-j}$ in this setting.

The operators $P$ and $P^\star$ are related in the following way. Denoting
by $U$ the symmetry of $B$ such that $U(b)_i = \sigma_i(b_{-i})$, we immediately
see that $U P^\star U = P$. We have
  $$P^\star P(b)_i= \frac{1}{(2d-1)^2} \sum_{ j\in I\atop{j\not= i}} \sum_{ k\in I\atop{k\not= j}} b_k$$
  and therefore
$$P^\star P = \frac{2d -2}{2d-1} UP + \frac{1}{2d-1} Id_B$$
By induction we get
$$(P^\star)^nP^n =  \frac{2d -2}{2d-1} UP^{2n-1} + \frac{1}{2d-1} (P^\star)^{n-1} P^{n-1}.$$
It follows that 
$$\aligned
P^{2n-1} &= \frac{2d -1}{2d-2} U(P^\star)^n P^n - \frac{1}{2d-2}U(P^\star)^{n-1} P^{n-1}\\
& =\frac{2d -1}{2d-2}  P^n (P^\star)^n U- \frac{1}{2d-2} P^{n-1} (P^\star)^{n-1}U.
\endaligned$$

Now, using corollary \ref{almunif}, we see that for every $b\in B$, the sequences
$\big(P^{2n-1}(b)\big)$ and $\big(P^{2n}(b)\big)$ converge almost uniformly. As a consequence  of the relation
(\ref{basrel}), we get, for $x\in A$, the almost uniform convergence of the sequences $\big(s_{2n-1}(x)\big)$,
 $\big(s_{2n}(x)\big)$ and therefore of $\big(\frac{1}{n}\sum_{k=0}^{n-1} s_k(x)\big)$.

Let $\E^{(2)}$ be the $\varphi$-invariant conditional expectation from $A$ onto the subalgebra $A^{(2)}$ 
of $\F_{d}^{(2)}$-invariant elements in $A$ and let us show that $\E^{(2)}(x)$ is the limit $\hat x$ of
 $\big(s_{2n}(x)\big)$. Since the automorphisms $\sigma_i$ are $\varphi$-preserving, we have
$\varphi(yx) = \varphi(ys_{2n}(x))$ for every $y \in A^{(2)}$. Moreover, the bounded sequence  $\big(s_{2n}(x)\big)$,
which converges to $\hat x$ almost uniformly, also converges to $\hat x$ in the
strong topology (see \cite[Theorem 1.1.3]{Jaj} for instance). Therefore we have, for $y\in A^{(2)}$,
$$\varphi(y\E^{(2)}(x)) =\varphi(yx) =\lim_{n\to +\infty}\varphi\big(ys_{2n}(x)\big) = \varphi(y\hat{x}).$$
To conclude that $\hat{x} = \E^{(2)}(x)$ it remains to check that $\hat x$ is
$\F_{d}^{(2)}$-invariant. We shall use the following identity:
\begin{equation}\label{libre}
s_1\circ s_n = \frac{2d-1}{2d} s_{n+1} + \frac{1}{2d}s_{n-1}.
\end{equation}
This relation is an immediate consequence of the fact that multiplying all words of length $1$
by all words of length $n$, each word of length $n+1$ is obtained once and each word of length $n-1$ is obtained $2d-1$ times.
It follows from (\ref{libre}) that
$$s_1^{2}\circ s_{2n} = \Big(\frac{2d-1}{2d}\Big)^2 s_{2n+2} +2\frac{2d-1}{(2d)^2}s_{2n} +
\frac{1}{(2d)^2} s_{2n-2},$$
from which we get $s_{1}^2(\hat x) = \hat x$.

If we view the elements of $A$ as vectors in the Hilbert space $L^2(A,\varphi)$ of the GNS representation associated
with $\varphi$, a strict convexity argument gives that $\sigma_w(\hat x) = \hat x$ for every word of length $2$.
Therefore, $\hat x$ is $\F_{d}^{(2)}$-invariant.

Similarly, we show that the limit of $\frac{1}{n}\sum_{k=0}^{n-1} s_k(x)$ is the conditional expectation
of $x$ with respect to the subalgebra of $\F_d$-invariant elements.
 
 The proof is similar for any $p>1$, still using corollary \ref{almunif}.
\end{proof}

\section{About the noncommutative Rota theorem}
 
 As said before, our corollary \ref{almunif} is a noncommutative version (for a particular Markov operator)
of the following theorem of G.-C. Rota \cite{Rot}~:
 
 \begin{thm}[\cite{Rot}\label{Rota}] Let $P$ be a $\mu$-preserving Markov operator on the pro\-bability measure space
 $(X,\mu)$. Then for
 every $p>1$ and $f\in L^p(X,\mu)$, the sequence $\big(P^n (P^\star)^n(f)\big)$ converges almost everywhere.
 \end{thm}
 
  This result was established by Rota from the
 commutative analogue of theorem  \ref{basic} (see \cite{Rot}). In fact,  Rota's theorem is still 
true whenever $f\in L Log L(X,\mu)$.

 The extension of Rota's theorem to the noncommutative setting (with a minimum of hypotheses on $P$)
is an interesting open problem.
   It can be studied using the general definition of $P^\star$ introduced in (\ref{adjoint1}). As a first step,
   we only consider the case of a state preserving Markov operator, so that $P^\star$
   satisfies (\ref{adjoint}).
   
 As  already explained in Section 4,  it is crucial to determine for which operators $P$ theorem \ref{basic} holds,
 and we shall examine hereafter this problem.

Let $(B,\phi)$ be a noncommutative probability space and $P: B \to B$ 
  be a $\phi$-preserving Markov operator. Obviously, if theorem
\ref{basic} holds, then $P$ must be factorizable in the following sense.

\begin{defn} \label{factor} Let $P : (A_1,\varphi_1) \to (A_0,\varphi_0)$ be a $(\varphi_0,\varphi_1)$-preserving
Markov operator. We say that
$P$ is {\it factorizable} if there exists a noncommutative probability space $(C, \psi)$ and two
normal unital  homomorphisms $j_0 : (A_0,\varphi_0) \to (C,\psi)$, 
$j_1: (A_1,\varphi_1) \to (C,\psi)$ such that
$(\psi,\varphi_0)$ and $(\psi,\varphi_1)$ are respectively $j_0$- and $j_1$-stationary 
and $P = j_{0}^\star\circ j_1$.
\end{defn}

Indeed, if theorem \ref{basic} holds, take $(C,\psi) = ({\mathcal B}, \Phi)$ and $j_i = J_i$, $i = 0,1$.
Since $J^\star_{0} = J_{0}^{-1} \circ \E_{0]}$, the equality $\E_{0]} \circ J_1 = J_0 \circ P$ gives
$P = j_{0}^{\star}\circ j_1$.

\begin{exs}\label{examples1}{\bf(a)} If $A_0$, $A_1$ are commutative abelian von Neumann algebras, every unital positive map
$P: A_1 \to A_0$ is factorizable. Indeed, let us endow the algebraic tensor product
$A_0\odot A_1$ with the inner product
$$\langle f_0\otimes f_1, g_0\otimes g_1 \rangle = \varphi_0\big(g_{0}^* f_0P(g_{1}^* f_1)\big)$$
and denote by $H_P$ the Hilbert space obtained by separation and completion.
Let $j_0$ (resp. $j_1$) be the normal representation of $A_0$ (resp.  $A_1$) defined by
$$j_0(g_0) (f_0\otimes f_1) = g_0f_0 \otimes f_1 \quad \hbox{for} \quad g_0, f_0 \in A_0, f_1 \in A_1$$
(resp.
$$j_1(g_1) (f_0\otimes f_1) = f_0 \otimes g_1 f_1 \quad \hbox{for} \quad  f_0 \in A_0, g_1, f_1 \in A_1).$$
We denote by $C$ the von Neumann algebra generated by $j_0(A_0)$ and $j_1(A_1)$ and by $\psi$
the state $c\mapsto \langle c \xi_P, \xi_P \rangle$ where $\xi_P$ is the class of $1\otimes 1$ in $H_P$.
In particular, we have $\psi\big(j_0(f_0)j_1(f_1)\big)= \varphi\big(f_0 P(f_1)\big)$. We easily check
that $P = j_{0}^\star\circ j_1$.

The reader familiar with Connes' notion of correspondence will reco\-gnize in $H_P$ the Hilbert space 
of the correspondence associated with $P$. This construction is also well known in probability
theory. In this framework, to $(A_i, \varphi_i)$ is associated a probability space $(X_i,\mu_i)$
such that $A_i = L^\infty(X_i,\mu_i)$ and $\dst\varphi_i(f) = \int_{X_i} f d\mu_i$ for $f\in A_i$,
$i=0,1$. Let $p$ be the probability transition associated with $P$, { \it i.e.} $\dst P(f)(x) = \int_{X_1} f(y) p(x,dy)$
for $f\in L^\infty(X_1,\mu_1)$. Then $C = L^\infty(X_0\times X_1,\nu)$ where $\nu$ is the
probability measure defined by
$$\int_{X_0\times X_1} f d\nu = \int_{X_0}\big(\int_{X_1} f(x_0,x_1)p(x_0, dx_1)\big) d\mu_0(x_0).$$
Here, $\psi$ is the state associated with $\nu$ and the embeddings $j_0$, $j_1$ are the obvious ones.

{\bf(b)} Let us come back to the setting $P: (B,\phi) \to (B,\phi)$ considered in section 4.
We take for $(C,\psi)$ the algebra $B_1 = \{b: I(2) \to A\}$ equipped with the state $\phi_1$.
Recall that $\dst\phi_1(b)= \sum_{w\in I(2)} p(w_0)p(w_0w_1)\varphi(b_{w_0 w_1})$ for
$b\in B_1$. Let $j_0$ be the map $\alpha_1 : B \to B_1$ defined by
$\alpha_1(b)_{w_0 w_1} = b_{w_0}$ and let $j_1$ be the map $\beta_1 : B \to B_1$ defined by
$\beta_1(b)_{w_0 w_1} = \sigma_{w_0}(b_{w_1})$. Then it is easily checked that $P = \alpha_{1}^\star\circ \beta_1$.

{\bf (c)} More generally, let us replace the endomorphisms $\sigma_i$ of $A$ by any factorizable 
$\varphi$-preserving completely positive map $P_i$ from $A$ to $A$ for $i\in I$ and let us
keep the notations of (b) for $(B,\phi)$. Then the map $P: (B,\phi) \to (B,\phi)$ defined
by 
$$P(b)_i = \sum_{j\in I} p(ij)P_i(b_j)$$
is still factorizable. Indeed, let us write $P_i = \alpha_{i}^\star \circ \beta_i$ for $\alpha_i, \beta_i$ 
appropriate maps from $(A,\varphi)$ into a noncommutative probability space $(C_i, \psi_i)$.
We define $C$ to be the von Neumann algebra formed of the
maps $b : w_0w_1 \in I(2) \mapsto b_{w_0 w_1} \in C_{w_0}$. We endow $C$ with the state
$$\psi : b \mapsto \sum_{w \in I(2)} p(w_0)p(w_0w_1) \psi_{w_0}(b_{w_0w_1}).$$
Moreover we define $\alpha, \beta$ from $B= A^I$ into $C$ by
$$\alpha(b)_{w_0w_1} = \alpha_{w_0}(b_{w_0}), \quad \beta(b)_{w_0w_1} = \beta_{w_0}(b_{w_1}).$$
Then we have $P = \alpha^\star\circ \beta$.

{\bf(d)} Another useful example is given by the following lemma.

\end{exs}

Before stating it, we make the convention hereafter, for simplicity of language, that all our maps between noncommutative
probability spaces  will (implicitely) be preserving Markov operators with respect to the given states.

\begin{lem} \label{libre1}  Let $(B,\varphi)$,  $(A_i,\varphi_i)$ $i= 0,1$, be non commutative probability spaces.
We assume that for $i=0,1$, $A_i$ contains a copy of $B$ as a von Neumann subalgebra 
and that there is a conditional expectation $\E_i : A_i  \to B$ with $\varphi_i = \varphi \circ \E_i$.
Then $P = i_0 \circ i_{1}^\star$ (where $i_k$, $k = 1,0$, is the canonical embedding of $B$ into $A_k$) is factorizable.
More precisely, one can find $(C,\psi)$ and homomorphisms $j_0$, $j_1$ from $(A_0,\varphi_0)$, $(A_1,\varphi_1)$
respectively, into $(C,\psi)$ such that $j_0\circ i_0 = j_1\circ i_1$ and $i_0 \circ i_{1}^\star = j_{0}^\star\circ j_1$.
\end{lem}

\begin{proof} We define $(C,\psi)$ to be the reduced amalgamated free product of $(A_1,\varphi_1)$
by $(A_0,\varphi_0)$ over $B$ and $j_0: A_0 \to C$, $j_1: A_1 \to C$ are the canonical
embeddings. This construction is explained for instance in \cite{BD}.  We give below
an equivalent construction, in terms of self-dual Hilbert modules (see \cite{Pas} for that notion).
It is the exact analogue of the construction exposed by Voiculescu in \cite[Section 5]{Voi}, except
that self-dual completions replace norm completions in the definitions ({\it e. g.} for relative tensor products).
For $i = 0,1$, we denote by $H_i$ the self-dual right Hilbert $B$-module obtained by completion of $A_i$ with respect to the 
B-valued inner product $\langle x, y\rangle = \E_i(x^* y)$ for $x,y \in A_i$ and $x \mapsto \hat{x}$
is the canonical embedding of $A_i$ into $H_i$. Recall
that left multiplication on
$A_i$ yields a unital normal injective homomorphism $\pi_i$ from $A_i$ into $\cB_B(H_i)$ 
(the von Neumann algebra of $B$-linear operators having an adjoint). In particular $H_i$ is a bimodule over $B$.
 We set $\xi_i = \widehat{1}_{A_i}$ and decompose $H_i$ as  the orthogonal direct sum of self-dual
Hilbert $B$ modules $H_i = \xi_i B\, \oplus \intH_i$. As in \cite{Voi} we define $(H,\xi) = (H_0, \xi_0)*(H_1,\xi_1)$
by
$$H = \xi B \,\oplus {\somtrois {n\geq 1} {\iota_1,\dots, \iota_n \in \{0,1\}}{\iota_1 \not= \iota_2
\not=\iota_3\not=\dots\not= \iota_{n-1}\not=\iota_n}}
 \intH_{\iota_1}\otimes_B \cdots \otimes_B \intH_{\iota_n}.$$
We denote by $j_i$ the normal representation of $A_i$ into $\cB_B(H)$ introduced by Voiculescu (written $\lambda_i$ in \cite{Voi})
and $C = A_0 *_B A_1$ is the von Neumann subalgebra of $\cB_B(H)$ generated by $j_0(A_0)$ and $j_1(A_1)$.
We shall freely identify the elements $a \in A_i$ and $j_i(a)\in C$, $i=0,1$. 

The map $\E : c \mapsto \langle \xi, c \xi\rangle$ is a normal conditional expectation from
$C$ onto $B$ and we set $\psi = \varphi \circ \E$. We identify the submodule $\xi B\, \oplus \intH_0$ of $H$ with $H_0$
and we denote by $q_0$ the orthogonal projection from $H$ onto $H_0$. Then $\widetilde{\E}_0 : c \mapsto
q_0 c q_0$ defines a normal conditional expectation from $C$ onto $j_0(A_0)$
  (via the identification $x\equiv j_0\circ \pi_{0}^{-1}(x)$ between $\pi_0(A_0)$ and $j_0(A_0)$)
such that $\widetilde{\E}_0 \circ j_1(a) = \E(a)$ for $a\in A_1$ (see \cite[Lemma 3.5]{BD}). 
 We have $\E = \E_0\circ \widetilde{\E}_0$ and therefore
$\varphi_0\circ  \widetilde{\E}_0 = \psi$.
Similarly we define $\widetilde{\E}_1 : C \to j_1(A_1)$.

It remains to show that $i_0\circ i_{1}^\star = j_{0}^\star\circ j_1$. For this, 
we shall compute $\langle \xi_2, j_{0}^\star\circ j_1(a_1)a_0 \xi_0 \rangle$ for $a_0\in A_0$ and $a_1\in A_1$.
We have, by a straightforward consequence of the definition of $j_0$ and $j_1$, 
\begin{align*}
\langle \xi_0, j_{0}^\star\circ j_1(a_1)a_0 \xi_0 \rangle & = \E_0\big( j_{0}^\star\circ j_1(a_1)a_0\big)\\
&= \E\big(j_1(a_1) j_0(a_0)\big)\\
&= \langle \xi, j_1(a_1) j_0(a_0)\xi \rangle\\
&= \langle j_1(a_1)^*\xi,  j_0(a_0)\xi \rangle\\
&= \langle \xi_1, a_1 \xi_1 \rangle \langle \xi_0, a_0 \xi_0 \rangle \\
&= \langle \xi_0, \langle \xi_1, a_1 \xi_1 \rangle  a_0 \xi_0 \rangle
\end{align*}
($\langle \xi_1, a_1\xi_1\rangle \in B$ being identified with $i_0\big(\langle \xi_1, a_1\xi_1\rangle\big)$ in the last
equality).
It follows that $j_{0}^\star\circ j_1(a_1) = i_0\circ i_{1}^\star(a_1)$ for every $a_1 \in A_1$.

\end{proof}

Note that $C$ is generated, as a von Neumann algebra, by $j_0(A_0)$ and $j_1(A_1)$.

\begin{thm}\label{basic1} A $\phi$-preserving Markov operator $P$ on $(B,\phi)$ is factorizable if and only 
if theorem \ref{basic} holds. As a consequence, the noncommutative Rota theorem \ref{almunif} applies to factorizable 
Markov operators.
\end{thm}

In fact, we shall prove a more complete theorem.

\begin{thm}[\textbf{Noncommutative Markov chain construction}]\label{ncMarkov} Let $P : (B,\phi) \to (B,\phi)$
be a factorizable $\phi$-preserving Markov operator. There exist
\begin{itemize}
\item [--] a von Neumann algebra $\mathcal{B}$,
 \item [--] a normal faithful state $\Phi$ on $\mathcal{B}$,
  \item [--] a normal unital  endomorphism $\beta : \mathcal{B} \to \mathcal{B}$
with $\Phi\circ \beta = \Phi$ and $\sigma_{t}^\Phi \circ \beta
= \beta \circ \sigma_{t}^\Phi$ for $t\in \R$,
 \item [--] a normal
unital homomorphism $J_0 : B \to \mathcal{B}$ with $\Phi\circ J_0 = \phi$ and $\sigma_{t}^\Phi \circ J_0
= J_0 \circ \sigma_{t}^\phi$ for $t\in \R$,
\end{itemize}
such that, if we set
  $J_n = \beta^n\circ J_0$ for $n\geq 0$, and if $\mathcal{B}_{n]}$ and $\mathcal{B}_{[n}$
denotes the von Neumann subalgebras of $\mathcal{B}$ generated by $\bigcup_{k\leq n}J_k(B)$
 and $\bigcup_{k\geq n}J_k(B)$ respectively, then
 \begin{itemize}
 \item[(i)] the algebras $\mathcal{B}_{n]}$ and $\mathcal{B}_{[n}$ are $\Phi$-invariant;
\item[(ii)] if $\mathbb{E}_{n]}$ and $\mathbb{E}_{[n}$ are the corresponding $\Phi$-preserving
conditional expectations,  we have
 \begin{align}\label{basicno1}
  &\mathbb{E}_{n]}\circ  J_q = J_n \circ P^{q-n},\quad \forall q\geq n\geq 0,\\
  &\mathbb{E}_{n +q]}\circ \beta^q = \beta^q \circ \E_{n]},\quad \forall n,q \in \N,\\
  &\mathbb{E}_{[n}\circ J_0 = J_n \circ (P^{\star})^n,\quad  \forall n\in \N.
\end{align}

 \end{itemize}
\end{thm}

\begin{proof} Assume the existence of $(B_1,\phi_1)$ and of homomorphisms
$\alpha_1$, $\beta_1$ from $(B,\phi)$ into $(B_1,\phi_1)$ with $P = \alpha_{1}^\star\circ \beta_1$. Let us
set $(B_0,\phi_0) = (B,\phi)$. Using lemma \ref{libre1}, 
we  construct inductively a sequence $(B_n,\phi_n)$ of noncommutative probability spaces,
and of pairs $(\alpha_{n+1},\beta_{n+1})$ of homomorphisms from $(B_n, \phi_n)$ into $(B_{n+1}, \phi_{n+1})$
such that 
\begin{align}\label{fondam}
\beta_n\circ \alpha_{n}^\star &= \alpha_{n+1}^\star\circ \beta_{n+1}\\
\alpha_{n+1}\circ \beta_{n} & = \beta_{n+1}\circ\alpha_{n}.
\end{align}
Let $({\mathcal B},\Phi)$ be the inductive limit of $(B_n,\phi_n)$ with respect to the morphisms
$\alpha_{n} : (B_{n-1},\phi_{n-1}) \to (B_n,\phi_n)$ and let $\lambda_n$ be the canonical embedding
from $(B_n,\phi_n)$ into $({\mathcal B},\Phi)$. We also define the endomorphism
$\beta$ of $({\mathcal B},\Phi)$, characterized by the condition
$$\lambda_{n+1} \circ \beta_{n+1} = \beta\circ \lambda_n\quad\hbox{for all}\quad n\geq 0.$$
Now we set $J_0 = \lambda_{0} : B \to {\mathcal B}$ and $J_n = \beta^n\circ J_0$ for $n\geq 0$.

We shall identify $B_n$ with the subalgebra $\lambda_n(B_n)$ of $\cB$ when convenient. Then, 
for $0\leq q\leq r$, the homomorphism $J_q$, when viewed as a homomorphism
from $B$ into $B_r \subset \cB$ is written
$$\aligned
J_q &= \alpha_r \circ\cdots\circ\alpha_{q+1}\circ\beta_q\circ\cdots\circ\beta_1\\
&= \beta_r \circ\cdots\circ\beta_{r-q+1}\circ\alpha_{r-q}\circ\cdots\circ\alpha_1,
\endaligned$$
if $q>0$. For $q = 0$,
$$J_0 = \alpha_r \circ\cdots\circ\alpha_1 .$$

We shall also need the following formula:
\begin{equation}\label{puissanceP}
P^n = \alpha_{1}^\star \circ\cdots\circ \alpha_{n}^\star\circ\beta_n\circ\cdots\circ\beta_1.
\end{equation}
Starting with $P= \alpha_{1}^\star\circ \beta_1$, this formula is proved by using (\ref{fondam}) repeatedly.

Since for every $n\geq 1$, the von Neumann algebra $B_n$ is generated by $\alpha_n(B_{n-1})$
and $\beta_n(B_{n-1})$, we see inductively that  
$$B_n \equiv \lambda_n(B_n) = B_{n]}.$$

Observe  that for $q\geq n\geq 0$, the conditional expectation $\E_{n]}$, restricted to $B_q$ is  
$\alpha_{n+1}^\star\circ \cdots \circ \alpha_{q}^\star$. It follows that 
$$\E_{n]} \circ J_q = \alpha_{n+1}^\star\circ \cdots \circ \alpha_{q}^\star\circ \beta_{q}\circ\cdots\circ\beta_1,$$
and, using (\ref{fondam}), we get
$$\aligned
\E_{n]} \circ J_q &= \beta_n \circ\cdots\circ\beta_1
\circ \alpha_{1}^\star\circ \cdots \circ \alpha_{q-n}^\star\circ \beta_{q-n}\circ\cdots\circ\beta_1\\
&= J_n\circ P^{q-n},
\endaligned$$
by (\ref{puissanceP}). This concludes the proof of (\ref{basicno1}).

Let us show now that $\E_{n+q]} \circ \beta^q = \beta^q \circ \E_{n]}$ for every $n,q\in \N$.
It is enough to fix an integer $r \geq n$ and to prove this formula in $B_r$.
Therefore we have to check that
$$\alpha_{n+q+1}^\star\circ \cdots \circ \alpha_{r+q}^\star\circ \beta_{r+q}\circ\cdots\circ\beta_{r+1}
= \beta_{n+q} \circ\cdots\circ \beta_{n+1}\circ\alpha_{n+1}^\star\circ\cdots\circ\alpha_{r}^\star,$$
which is still a consequence of (\ref{fondam}).

Finally, let us prove  that $\E_{[n} \circ J_0 = J_n \circ (P^\star)^n$. We have to check that 
for $b\in B$, $r\geq n$ and $ c\in \Big( \bigcup_{n\leq k \leq r} J_k(B)\Big)''$ then
$$\Phi\big(cJ_0(b)\big) = \Phi \big(cJ_n \circ (P^\star)^n(b)\big).$$
In order to shorten notations, for any pair of integers $k\leq \ell$, we shall write
$\alpha_{[\ell,k]}$ for $\alpha_\ell\circ\cdots\circ\alpha_k$, and we shall adopt the same
convention for $\beta_{[\ell,k]}$. We denote by $B_{[n,r]}$ the algebra $\Big( \bigcup_{n\leq k \leq r} J_k(B)\Big)''$
when viewed into $B_r$. 
One shows immediately that for every $1\leq k \leq r$,
$$B_{[r-k,r]} =\beta_{[r,k+1]}(B_k).$$

We have to prove that for $n\leq r$, 
\begin{equation}\label{induct}\E^{B_r}_{\beta_{[r,r-n+1]}(B_{r-n})}\circ \alpha_{[r,1]} = 
\alpha_{[r,n+1]}\circ \beta_{[n,1]}\circ (P^\star)^n,
\end{equation}
where $\E^{B_r}_{\beta_{[r,r-n+1]}(B_{r-n})}$ is the natural conditional expectation from $B_r$
onto $\beta_{[r,r-n+1]}(B_{r-n}) = B_{[n,r]}$.
We fix $n$ and proceed by induction on $r$. 

First, for $r=n$, let us prove that $\E^{B_n}_{\beta_{[n,1]}(B)}\circ\alpha_{[n,1]} =  \beta_{[n,1]}\circ (P^\star)^n$:
by (\ref{puissanceP}) we have $(P^\star)^n= \beta_{[n,1]}^\star\circ \alpha_{[n,1]}$
and therefore
$$ \beta_{[n,1]}\circ (P^\star)^n =\beta_{[n,1]}\circ  \beta_{[n,1]}^\star\circ \alpha_{[n,1]} = \E^{B_n}_{\beta_{[n,1]}(B)}\circ \alpha_{[n,1]}.$$

Asssume now that (\ref{induct}) is true and let us prove the formula for $r+1$. Take
$c\in B_{r+1-n}$ and $b\in B$ and compute 
$$\phi_{r+1}\big(\beta_{[r+1, r-n+2]}(c) \alpha_{[r+1,1]}(b)\big) = \phi_r\big(\alpha_{r+1}^\star\circ \beta_{[r+1, r-n+2]}(c) \alpha_{[r,1]}(b)\big).$$ 
Using again (\ref{fondam}), we get
$$ \alpha_{r+1}^\star\circ \beta_{[r+1, r-n+2]}(c) =\beta_{[r,r-n+1]}\circ \alpha_{r-n+1}^\star(c).$$
We have $\alpha_{r-n+1}^\star(c)\in B_{r-n}$ and therefore the assumption (\ref{induct}) gives
\begin{align*}
\phi_{r+1}\big(&\beta_{[r+1, r-n+2]}(c)\alpha_{[r+1,1]}(b)\big)\\
 &=\phi_r\big(\alpha_{r+1}^\star\circ \beta_{[r+1, r-n+2]}(c)\alpha_{[r,1]}(b)\big)  \\
&= \phi_r\big(\alpha_{r+1}^\star\circ \beta_{[r+1, r-n+2]}(c)\E^{B_r}_{\beta_{[r,r-n+1]}(B_{r-n})}\circ \alpha_{[r,1]}(b)\big)\\
 & = \phi_r\big(\alpha_{r+1}^\star\circ\beta_{[r+1, r-n+2]}(c) \alpha_{[r,n+1]}\circ \beta_{[n,1]}\circ (P^\star)^n(b)\big)\\
 &= \phi_{r+1}\big(\beta_{[r+1, r-n+2]}(c)\alpha_{r+1}\circ \alpha_{[r,n+1]}\circ\beta_{[n,1]}\circ (P^\star)^n(b)\big).
 \end{align*}
 It follows that 
 $$\E^{B_{r+1}}_{\beta_{[r+1,r-n+2]}(B_{r+1-n}})\circ  \alpha_{[r+1,1]}(b) = \alpha_{[r+1,n+1]}\circ \beta_{[n,1]}\circ
(P^\star)^n(b).$$ This concludes the proof of theorem \ref{basic1}.
\end{proof}

\begin{rem} {\bf (a)} Let us consider the commutative case (see (\ref{examples1}) (a)).
Given a $\mu$-preserving Markov operator $P$ from $L^\infty(X,\mu)$ into itself, there is a much simpler
construction of $(\cB, \Phi, \beta, J_n)$ satisfying theorem \ref{ncMarkov}, known as the Daniell-Kolmogorov
construction. The von Neumann $\cB$ is $L^\infty(\Omega,\nu)$, where $\Omega = X^\n$ is the trajectory space
equipped with the usual $\sigma$-field of measurable subsets, and $\nu$ is the Markov measure associated
with $P$ and the starting measure $\mu$. Let $(X_n)_{n\in \N}$ be the corresponding
homogeneous Markov chain, that is, for $n\in \N$ the random variable $X_n : \Omega \to X$ is the projection
$(x_i)_{i\in \n} \mapsto x_n$. Then $J_n$ is the homomorphism $f \mapsto f\circ X_n$ from $L^\infty(X,\mu)$
into $L^\infty(\Omega, \nu)$.

{\bf (b)} When $A$ is the von Neumann algebra $\C$, the construction made in section 4 is a particular case
of the above one. It corresponds to the case where $X = I$, $\mu$ is the distribution $\big(p(i)\big)_{i\in I}$
and $P : \C^I \to \C^I$ is given by the transition matrix $[p(ij)]$.

{\bf (c)} Given any unital completely positive map $P : B \to B$, where $B$ is a von Neumann algebra
(or a unital $C^\star$-algebra), without any other assumption on $P$, several Daniell-Kolmogorov type constructions
can be found in the literature (see \cite{Sau1}, \cite{BP}).  In these constructions, equalities (\ref{basicno1}) and
(6.2) are fulfilled. The main point here is that we also need condition (6.3).
\end{rem}

\noindent{\bf Problems.} (1) Are there examples of $\phi$-preserving Markov operators not being factorizable?

(2) For a factorizable Markov operator $P$, find a canonical construction for $(\cB, \Phi,\beta, J_n)$ in theorem
\ref{ncMarkov}.

(3) Let $P$ be a normal unital completely positive map from $(B,\phi)$ onto itself such that $\phi\circ P = \phi$.
Let $P^\star$ be the completely positive operator such that
$\phi\big(P^\star(b) \sigma_{-i/2}^\phi(a)\big)= \phi\big(\sigma_{i/2}^\phi(b)P(a) \big)$
for all analytic elements $a,b\in B$. Study the convergence of the sequence $\big(P^n(P^\star)^n\big)$.

\begin{rem} Junge and Xu have generalized in \cite{JX3} both the individual Dun\-ford-Schwartz  ergodic theorem 
and the classical Stein ergodic theorem proved in \cite{Ste}. Using our
terminology, the noncommutative Stein ergodic theorem reads as follows:

\begin{thm}[\textbf{Junge-Xu}, Corollary 7.11, \cite{JX3}]
 Let $(B,\phi)$ be a noncommutative probability space and let $P: B\to B$
be a $\phi$-preserving Markov operator such that $P = P^\star$. Then for $p\in ]1,+\infty]$ and
$x\in L^p(B,\phi)$, the sequence $\big(P^n(x)\big)$  converges bilaterally almost surely to the canonical
projection of $x$ onto the subspace of $P$-invariant elements. If $p>2$, the convergence
also holds almost surely.
\end{thm}

This gives of course the noncommutative Rota theorem for self-adjoint $\phi$-preserving
Markov operators.
\end{rem}

\bibliographystyle{plain}

 \end{document}